\documentclass[graybox]{svmult}

\usepackage{amsmath,amssymb,bbm}
\usepackage{setspace} 

\usepackage{amsfonts}
\usepackage{mathptmx}
\usepackage{mathptmx}       
\usepackage{helvet}         
\usepackage{courier}        
\usepackage{type1cm}        
\usepackage{makeidx}         
\usepackage{graphicx}        
\usepackage{multicol}        
\usepackage[bottom]{footmisc}

\makeindex             

\begin{document}

\def \\ { \cr }
\def\R{\mathbb{R}}
\def\N{\mathbb{N}}
\def \d{{\rm d}}

\def\QED{\hfill $\Box$}

\title*{Local times for functions with finite variation:  two versions of Stieltjes change of variables formula}
\titlerunning{Local times for functions with finite variation} 
\author{Jean Bertoin and Marc Yor}
\institute{Jean Bertoin \at Institut f\"ur Mathematik, 
Universit\"at Z\"urich, 
Winterthurerstrasse 190, 
CH-8057 Z\"urich, Switzerland, \email{jean.bertoin@math.uzh.ch}
\and Marc Yor \at Institut Universitaire de France and Laboratoire de Probabilit\'es et Mod\`eles Al\'eatoires, 
UPMC, 4 Place Jussieu, 75252 Paris cedex 05, France \email{deaproba@proba.jussieu.fr}}
%
%
\maketitle

\abstract { We introduce two natural notions for the occupation measure of a function $V$ with finite variation. The first yields a signed measure, and the second a positive measure. By comparing two versions of the change-of-variables formula, we show that both measures are absolutely continuous with respect to Lebesgue measure. Occupation densities can be thought of as local times of $V$, and are described by a Meyer-Tanaka like formula.}



\hskip 4mm

\section{Introduction}

Although the setting of the present work is entirely deterministic and does not involve Probability Theory at all,  its motivation comes from an important chapter of Stochastic Calculus
(see, for instance, Section IV.7 in Protter \cite{Pr}).
Specifically, let $(X_t)_{t\geq 0}$ be a real-valued semimartingale, i.e. $X=M+V$ where $M$ is a local martingale and $V$ a process with finite variation, and the paths of $X$ are right-continuous and possess left-limits (c\`adl\`ag) a.s. In turn the local martingale  can be expressed as the sum $M=M^c+ M^d$ of a continuous local martingale $M^c$ and a purely discontinuous local martingale $M^d$. The occupation measure of $X$ on a time interval $[0,t]$
is defined as
\begin{equation} \label{OMm}
\int_0^t {\bf 1}_A(X_s) \d \langle M^c\rangle_s\,,\qquad A\in{\mathcal B}(\R)\,,
\end{equation}
where $\langle M^c\rangle$ denotes the quadratic variation of $M^c$. 
A fundamental result in this area is that this occupation measure is a.s. absolutely continuous with respect to Lebesgue measure, i.e. the following occupation density formula holds:
$$\int_0^t {\bf 1}_A(X_s) \d \langle M^c\rangle_s = \int_A L^x_t\d x\,. $$
The occupation densities $\{L^x_t: x\in \R \hbox{ and } t\geq 0\}$ are known as the local times of $X$ and are described by the Meyer-Tanaka formula
\begin{eqnarray}\label{E0}
&&(X_t-x)^+-(X_0-x)^+  \nonumber \\
&=& \int_0^t {\bf 1}_{\{X_{s-}>x\}} \d X_s  + \frac{1}{2} L^x_t 
+ \sum_{0< s \leq t}\left( {\bf 1}_{\{X_{s-}\leq x\}} (X_s-x)^+ + {\bf 1}_{\{X_{s-}>x\}}(X_s-x)^-\right).
\end{eqnarray}

The purpose of this work is to point out that the above results  for semimartingales  have natural analogs in the deterministic world of functions with finite variation. 
The rest of this note is organized as follows. Our main result is stated in Section 2 and proven in Section 4. Section 3 is devoted to two versions of the Stieltjes change-of-variables formula which lie at the heart of the analysis.

\section{Local times for functions with finite variation}
We consider a c\`adl\`ag function $V: [0,\infty)\to \R$ with finite variation. There is  the canonical decomposition of $V$ as the sum of its continuous and discontinuous components,
$$V= V^{c}+V^{d},$$ where 
$$V^{d} (t) =\sum_{0< s \leq t}Ê\Delta V(s)$$
and the series is absolutely convergent. 
Plainly, the continuous part $V^c$ has also finite variation;  we write $V^c(\d t)$ for the corresponding signed Stieltjes measure and  $|V^c(\d t)|$ for its total variation measure. 

We then  introduce for every $t>0$  a signed measure $ \theta_t(\d x)$ and a positive measure $\vartheta_t(\d x)$, which are defined respectively by
\begin{equation}\label{OCV} 
\theta_t(A) = \int_0^t {\bf 1}_A (V(s)) V^c(\d s)
\end{equation}
and
\begin{equation}\label{OCA} 
 \vartheta_t(A) = \int_0^t {\bf 1}_A (V(s)) |V^c(\d s)|
 \end{equation}
where $A\subseteq \R$ stands for a generic Borel set. We interpret $ \theta_t(\d x)$ and $\vartheta_t(\d x)$ respectively as the signed and the absolute occupation measure of $V$ on the time interval $[0,t]$. 

Our goal is to show that the occupation measures $\theta_t(\d x)$ and $\vartheta_t(\d x)$ are both absolutely continuous with respect to the Lebesgue measure
and to describe  explicitly their densities. In this direction, we introduce the following definitions and notations. 

For every $x\in\R$ and $t>0$, we say that $V$ increases through the level $x$ at time $t$ and then write $t\in{\mathcal I}(x)$ if :
\begin{itemize}
\item $V(t)=x$ and  $V$ is continuous at time $t$, 
\item $V(s)-V(t)$ has the same sign as $s-t$ for all $s$ in some neighborhood of $t$. 
\end{itemize}
Similarly, we say that $V$ decreases through the level $x$ at time $t$ and then write $t\in{\mathcal D}(x)$ if $-V$ increases through the level $-x$ at time $t$. We then define 
$$\ell^x(t)={\rm Card}\{(0,t]\cap {\mathcal I}(x)\}-{\rm Card}\{(0,t]\cap {\mathcal D}(x)\}\,,$$
whenever the two quantities in the difference in the right-hand side are finite.
For each such $x$, $t\to \ell^x(t)$ is a right-continuous function with integer values, which can only jump at times when $V$ reaches $x$ continuously.
 We also set
$$\lambda^x(t)={\rm Card}\{(0,t]\cap {\mathcal I}(x)\}+{\rm Card}\{(0,t]\cap {\mathcal D}(x)\}\,,$$
so that the point measure $\lambda^x(\d t)$ coincides with the total variation measure $|\ell^x(\d t)|$.

Our main result identifies  $(\ell^x(t): x\in\R \hbox{ and } t\geq 0)$ and $(\lambda^x(t): x\in\R \hbox{ and } t\geq 0)$ respectively as the signed and absolute local times of $V$. 
\begin{theorem} \label{T1} For every $t>0$, 
there are the identities
$$\theta_t(\d x) = \ell^x(t) \d x \quad \hbox{and} \quad \vartheta_t(\d x) = \lambda^x(t) \d x.$$
\end{theorem} 
 
Theorem \ref{T1} is an immediate consequence of the chain rule when $V$ is piecewise of class ${\mathcal C}^1$ (with possible jumps); however the general case requires a more delicate analysis. 
It is of course implicit in the statement that $\ell^x(t)$ and $ \lambda^x(t)$ are well-defined for almost all $x\in\R$ and yield measurable functions in $L^1(\d x)$, which is not {\it a priori} obvious.  

It is  interesting to point out that essentially the same result was established by Fitzsimmons and Port \cite{FP} in the special case when $V$ is a L\'evy process with finite variation and non-zero drift (then, when for instance the drift is positive, only increasing passages through a level may occur). 
In the same direction, Theorem \ref{T1} applies much more generally  to sample paths of a large class of (semi) martingales with finite variation, which arise for instance by predictable compensation of a pure jump process with finite variation. Note that when $M$ is a local martingale with finite variation, then its continuous component in the sense of Martingale Theory is always degenerate (i.e. identically zero), whereas its continuous component in the sense of functions with finite variation may be non-trivial. In that case, the occupation measure defined by \eqref{OMm} is degenerate, and then definitions \eqref{OCV} and \eqref{OCA} may be sounder. 

We further stress that for semi-martingales with non-degenerate continuous martingale component $M^c\not \equiv 0$,  the  local time $t\to L^x_t$ defines an increasing process which is always  continuous, whereas for functions with finite variation, $t\to \ell^x(t)$ and $t\to \lambda^x(t)$ are integer-valued step-functions which only jump at times when the function $V$ crosses continuously the level $x$. 

We also observe the following consequence of Theorem \ref{T1}. Provided that  $V$ is not purely discontinuous, i.e. $V^c\not \equiv 0$, then the set of levels through which $V$ increases or decreases has a positive Lebesgue measure. This contrasts sharply with the case of a typical Brownian path, or also, say, a sample path of a symmetric stable L\'evy process with index $\alpha\in(1,2)$,  as then local times exist, but there is no level through which the sample path increases or decreases (cf. Dvoretsky {\it et al.} \cite{DEK} and \cite{Be}). Of course, in those examples, sample paths have infinite variation. 

\section{Two versions of the change-of-variables formula}
Our strategy for proving Theorem \ref{T1} is to compare two change-of-variables formulas. The first version is standard (see, for instance, Dellacherie and Meyer \cite{DM}  on pages 168-171); a proof will be given for the sake of completeness.

\begin{proposition} \label{P1} 
Let $f:\R\to \R$ be a function of class ${\mathcal C}^1$. Then for every $t>0$,  there is the identity
$$f(V(t))-f(V(0))=\int_0^t f'(V(s)) V^c(\d s) + \sum_{0< s \leq t } \left(f(V(s))-f(V(s-))\right)\,.$$
\end{proposition}
\proof
With no loss of generality, we may assume that $f'$ has compact support with $|f'(x)|\leq 1$ for all $x\in\R$. 
We fix $\varepsilon >0$ arbitrarily small and may find a finite sequence of times $0=t_0< t_1< \cdots < t_n = t$ such that 
\begin{equation}\label{E1}
\max_{i=1, \ldots , n}\  \max_{t_{i-1}\leq u,v \leq t_i}|V^c(u)-V^c(v)| \leq \varepsilon
\end{equation}
 and
\begin{equation}\label{E2}
\sum_{\substack{0< s <t\\s\neq t_1, \ldots, t_n}}  |\Delta V(s)| \leq \varepsilon\,
\end{equation}
where, as usual, $\Delta V(s)=V(s)-V(s-)$. 
We write 
$$f(V(t))-f(V(0))= \sum_{i=1}^n \left(f(V(t_i-))-f(V(t_{i-1}))\right) + \sum_{i=1}^n \left(f(V(t_{i}))-f(V(t_{i}-))\right).$$
On the one hand, by the triangle inequality, we have 
\begin{eqnarray*} \left| \sum_{\substack{0< s <t\\s\neq t_1, \ldots, t_n}} \left(f(V(s))-f(V(s-))\right)\right|
&\leq &\sum_{\substack{0< s <t\\s\neq t_1, \ldots, t_n}}  \left| f(V(s))-f(V(s-))\right|\\
&\leq &\sum_{\substack{0< s <t\\s\neq t_1, \ldots, t_n}}  |\Delta V(s)|\,,
\end{eqnarray*}
where in the last line, we applied the mean value theorem and the assumption that $|f'|\leq 1$. Using \eqref{E2},  we have therefore shown that $$ \left| \sum_{0< s \leq t } \left(f(V(s))-f(V(s-))\right)- \sum_{i=1}^n \left(f(V(t_{i}))-f(V(t_{i}-))\right)\right| \leq \varepsilon.$$

On the other hand, again from the mean value theorem,  for each $i=1, \ldots , n$, there is some $x_i$ between $V(t_{i-1})$ and $V(t_{i}-)$ such that 
$$f(V(t_i-))-f(V(t_{i-1}))= f'(x_i)\left(V^c(t_i)-V^c(t_{i-1})\right) +  f'(x_i)\left(V^d(t_i-)-V^d(t_{i-1})\right).$$
Since $|f'|\leq 1$, it follows readily from \eqref{E2} and the triangle inequality  that
 $$ \sum_{i=1}^n\left|f'(x_i)(V^d(t_i-)-V^d(t_{i-1}))\right| \leq \varepsilon.$$
 Further, since $x_i$ lies between $V(t_{i-1})$ and $V(t_{i}-)$, we deduce from \eqref{E1} and \eqref{E2} that
$$|x_i-V(s)|\leq 2\varepsilon\quad \hbox{ for every }s\in[t_{i-1}, t_i).$$
 It follows that  
 $$\left|  \sum_{i=1}^nf'(x_i)(V^c(t_i)-V^c(t_{i-1})) - \int_0^t f'(V(s))V^c(\d s)Ê\right| \leq W^c(t) \sup_{|x-y|\leq 2\varepsilon}|f'(x)-f'(y)| ,
 $$
 where $W^c(t)<\infty$ stands for the total variation of $V^c$ on $[0,t]$. 

Putting the pieces together, we arrive at the inequality
\begin{eqnarray*}
&&\left| f(V(t))-f(V(0))-\int_0^t f'(V(s)) V^c(\d s) - \sum_{0< s \leq t } \left(f(V(s))-f(V(s-))\right)\right|\\
&\leq& 2\varepsilon +W^c(t)  \sup_{|x-y|\leq 2\varepsilon}|f'(x)-f'(y)| .
\end{eqnarray*}
Because $f'$ is continuous and has compact support, the upperbound tends to $0$ when $\varepsilon\to 0+$, which establishes the change-of-variables formula. \QED

In order to state a second version of the change-of-variables formula, we need first to introduce some further notions.
We say that $x\in\R$ is a {\it simple level} for $V$ if the set 
$$\{t>0: V(t-)<x<V(t) \hbox{ or } V(t)< x <V(t-) \hbox{ or } V(t)=x \}$$ 
is  discrete and if there is no jump of $V$ that starts or ends at $x$, i.e.
$$\{t>0: \Delta V(t)\neq 0 \hbox{ and either }x=V(t) \hbox{ or } x=V(t-)\}=\varnothing.$$
Otherwise, we say that $x$ is a {\it complex level} for $V$.

\begin{lemma} \label{L1} 
The set of complex levels for $V$ has zero Lebesgue measure.
\end{lemma}

\proof 
The set of jump times of $V$ is at most countable. Thus so is 
$$\{x\in\R: x=V(t) \hbox{ or } x=V(t-) \hbox{ for some $t>0$ with } \Delta V(t) \neq 0\}\,,$$
and {\it a fortiori} this set has zero Lebesgue measure.

Consider a sequence of partitions $\tau_n$ of $[0,\infty)$ with mesh tending to $0$ and such that $\tau_{n+1}$ is thiner as $\tau_n$. We write $0=t_0^n< t_1^n < \ldots$ for the elements of $\tau_n$. We set for every $n, i\in\N$
 $$J_i^n=[\inf_{t_i^n\leq s \leq t_{i+1}^n}V(s), \sup_{t_i^n\leq s \leq t_{i+1}^n}V(s)]\,,$$
 and write as usual $|J_i^n|$ for the length of this interval. If we write $W(t)$ for the total variation of $V$ on $[0,t]$, then   for every  $t\in\tau_n$, there is the upperbound
 $$\sum_{i: t_i^n< t}|J_i^n| \leq W(t). $$
 
 Now introduce 
 $$N_x(t) = {\rm Card} \{0< s \leq t: V(s-)\leq x \leq V(s) \hbox{ or } V(s)\leq x \leq V(s-)\}$$
 and observe that 
 $$N_x(t) = \lim_{n\to \infty}{\uparrow} \,   \sum_{i: t_i^n <  t} {\bf 1}_{ J_i^n}(x). $$
It follows from monotone convergence and Fubini-Tonelli  Theorem that
\begin{equation}\label{E3}
\int_{\R} N_x(t) \d x \leq W(t) <\infty.
\end{equation} 
{\it A fortiori} $N_x(t)<\infty$ for almost every $x$, which entails that the set 
 $$\{t>0: V(t-)\leq x \leq V(t) \hbox{ or } V(t)\leq x \leq V(t-)\}$$ is discrete for almost every $x$.  \QED
 
  We stress that for every simple level $x$, we can enumerate the increasing and the decreasing passage times of $V$ through $x$, and $\ell^x(t)$ and $\lambda^x(t)$ are both well-defined and finite.  We are now able to state the second version of the change-of-variables formula. 
 
\begin{proposition} \label{P2} 
Let $f:\R\to \R$ be a function of class ${\mathcal C}^1$. Then for every $t>0$,  there is the identity
$$f(V(t))-f(V(0))=\int_{\R} f'(x) \ell^x(t) \d x + \sum_{0< s \leq t } \left(f(V(s))-f(V(s-))\right)\,.$$
\end{proposition}
{\bf Remark.} Our  recent note \cite{BY} (see Theorem 2.1(i) there) contains a similar result when $V$ is an increasing process and the function $f$ is merely assumed to be non-decreasing. Of course, when $V$ is non-decreasing, the number of increasing passages through any level is at most $1$ and there are no decreasing passages, so $\ell^x(t) = 1$ or $0$, depending essentially on whether $V$ hits the level $x$ before time $t$ or not.
\proof
Consider a simple level $x$. We further suppose that $ x\neq V(t)$ and $x\neq V(0)$. 
For any such $x$, one can define the consecutive crossings of the level $x$ by the graph of $V$.
Each crossing can be made either upwards or downwards, and occurs either continuously or by a jump.
More precisely, an upward crossing of $x$ corresponds either to an increase through the level $x$, or to a positive jump starting strictly below $x$ and finishing strictly above $x$; and there is an analogous alternative for downwards crossings. Plainly upwards and downwards crossings alternate,
and the quantity
$${\bf 1}_{[x,\infty)}V(t) - {\bf 1}_{[x,\infty)}V(0)$$
coincides with the difference between upwards and downwards crossings of the level $x$ made before time $t$. Distinguishing crossings according to whether they occur continuously or by a jump, we arrive at
the identity
\begin{equation} \label{E4}
{\bf 1}_{[x,\infty)}V(t) = {\bf 1}_{[x,\infty)}V(0)+ \ell^x(t) + \sum_{0< s \leq t} \left({\bf 1}_{[x,\infty)}V(s)- {\bf 1}_{[x,\infty)}V(s-)\right)\,,
\end{equation}
where the series in the right-hand side accounts for the difference between the number of positive jumps and the number of negative jumps  across the level  $x$. 

It is convenient to write $\nu$ for the point measure on $[0,\infty)$ which has an atom at each jump time of $V$, viz.
$$\nu(\d t) = \sum_{s\geq 0} {\bf 1}_{\{\Delta V(s)\neq 0\}} \delta_s(\d t)\,.$$
We stress that $\nu$ is sigma-finite (because the number of jump times of $V$ is at most countable), and rewrite 
$$\sum_{0< s \leq t} \left({\bf 1}_{[x,\infty)}V(s)- {\bf 1}_{[x,\infty)}V(s-)\right)=\int_{[0,t]} \left({\bf 1}_{[x,\infty)}V(s)- {\bf 1}_{[x,\infty)}V(s-)\right)\nu(\d s).$$
We also observe that, in the notation of the proof of Lemma \ref{L1}, there is the inequality
$$\int_{[0,t]} \left|{\bf 1}_{[x,\infty)}V(s)- {\bf 1}_{[x,\infty)}V(s-)\right|\nu(\d s) \leq N_x(t)$$
and recall from \eqref{E3} that the function $x\to N_x(t)$ is in $L^1(\d x)$. Because the map 
$$(s,x)\to \left({\bf 1}_{[x,\infty)}V(s)- {\bf 1}_{[x,\infty)}V(s-)\right)$$
is measurable, it follows from Fubini Theorem that the map
$$x\to \int_{[0,t]} \left({\bf 1}_{[x,\infty)}V(s)- {\bf 1}_{[x,\infty)}V(s-)\right)\nu(\d s)$$
is also measurable. In particular, we conclude  from \eqref{E4} that $x\to \ell^x(t)$ is measurable. 

We now rewrite the identity  \eqref{E4} in the form
$$
{\bf 1}_{[x,\infty)}V(t) = {\bf 1}_{[x,\infty)}V(0)+ \ell^x(t) + \int_{[0,t]} \left({\bf 1}_{[x,\infty)}V(s)- {\bf 1}_{[x,\infty)}V(s-)\right)\nu(\d s)\,,
$$
 multiply both sides  by $f'(x)$ and integrate with respect to the Lebesgue measure $\d x$. An application of Fubini Theorem (which is legitimate thanks to the observations above) yields
$$
f(V(t))=f(V(0))+ \int_{\R} f'(x) \ell^x(t) \d x + \int_{[0,t]}\left(f(V(s))-f(V(s-))\right) \nu(\d s)\,,$$
which is the change-of-variables formula of the statement. \QED 

\section{Proof of Theorem \ref{T1}}

We now turn our attention to  the proof of our main result, which will also require the following elementary observation.

\begin{lemma} \label{L2}
Let $\mu$ be a signed measure on some measurable space, and $\varphi$ a measurable function with values in $\{-1, 1\}$ such that $\varphi \mu$ is a positive measure. Then $\varphi \mu$ coincides with the total variation measure $|\mu|$.
\end{lemma}
\proof Indeed, $\varphi \mu = |\varphi \mu|$ since $\varphi \mu$ is a positive measure, and on the other hand, we have also  $|\varphi \mu|=|\varphi| | \mu| =|\mu|$ since $ |\varphi|\equiv 1$. \QED 

We now have all the ingredients needed for establishing Theorem \ref{T1}. 

{\noindent{\bf Proof of Theorem \ref{T1} :}\hskip10pt}  
First, comparing the two change-of-variables formulas in Propositions \ref{P1} and \ref{P2}, we get that for any continuous function $g:\R\to \R$
(think of $g=f'$ as the derivative of a ${\mathcal C}^1$ function),
there is the identity
$$\int_0^t g(V(s)) V^c(\d s) = \int_{\R} g(x) \theta_t(\d x) = \int_{\R} g(x) \ell^x(t) \d x.$$
Thus $\theta_t(\d x)= \ell^x(t) \d x$. 

Next, we introduce the signed measure $\mu(\d s, \d x)$ on $[0,t]\times \R$ which is defined  by
$$\mu(A)= \int_0^t  {\bf 1}_A (s,V(s))  V^c(\d s)\,, \qquad A\in{\mathcal B}([0,t]\times \R).$$
The (signed) occupation density formula above  yields  that for every $f:\R\to \R$ and $g: [0,t]\to \R$ measurable and bounded, there is the identity
$$\int_{[0,t]\times \R} f(x) g(s) \mu(\d s, \d x) = \int_0^t f(V(s)) g(s) V^c(\d s) = \int_{\R} \left( \int_{[0,t]} g(s) \ell^x(\d s)\right)  f(x) \d x.$$
More precisely,  when $g$ is a step function, this identity follows from the occupation density formula  by linearity, and the general case is then derived through a version of the  monotone class theorem (cf. Neveu \cite{Ne}).

Then consider a measurable function $\varphi: [0,\infty) \to \{1,-1\}$ such that $V^c(\d t ) = \varphi(t) |V^c(\d t)|$. 
We have 
$$\int_{[0,t]\times \R} f(x) g(s) \varphi(s)  \mu(\d s, \d x) = \int_0^t f(V(s)) g(s) |V^c(\d s)|,$$
and the right-hand side is nonnegative whenever $f,g\geq 0$. Again by a version of the monotone class theorem, this shows that $ \varphi(s)  \mu(\d s, \d x)$ is a positive measure, and more precisely, since $|\varphi|=1$, Lemma \ref{L2} shows that $ \varphi(s)  \mu(\d s, \d x)= | \mu(\d s, \d x)|$ is the total variation measure of $ \mu(\d s, \d x)$. Now we write
$$\int_0^t f(V(s)) g(s) |V^c(\d s)| = \int_{\R} \left( \int_{[0,t]} g(s) \varphi(s)  \ell^x(\d s)\right)  f(x) \d x.$$
Plainly, whenever $g\geq 0$, we must have $\int_{[0,t]} g(s) \varphi(s)  \ell^x(\d s)\geq 0$ for almost all $x\in \R$, that is $\varphi(s)  \ell^x(\d s)$ is a positive measure for almost all $x\in \R$. Again because $|\varphi|=1$, this entails that 
$$\varphi(s)  \ell^x(\d s) = | \ell^x(\d s)| = \lambda^x(\d s)\,,$$
where the second equality is merely the definition of $\lambda^x$. Putting the pieces together, we have shown that
 $$\int_0^t f(V(s)) g(s) |V^c(\d s)| = \int_{\R} \left( \int_{[0,t]} g(s) \lambda^x(\d s)\right)  f(x) \d x,$$
 which for $g\equiv 1$ simply reads $\vartheta_t(\d x) = \lambda^x(t) \d x$. \QED
 
 We now conclude this note by stressing that the elementary identity \eqref{E4} for the signed local time $\ell^x(t)$ at a simple level $x$ should be viewed as the analog of the Meyer-Tanaka formula  \eqref{E0} in the semimartingale setting. We also point at the alternative formula
 \begin{equation}\label{TK1}
 {\bf 1}_{(-\infty,x)}V(t) = {\bf 1}_{(-\infty,x)}V(0)- \ell^x(t) + \sum_{0< s \leq t} \left({\bf 1}_{(-\infty,x)}V(s)- {\bf 1}_{(-\infty,x)}V(s-)\right)\,.
 \end{equation}

 For the absolute local time $\lambda^x(t)$, one sees similarly that for every simple level $x$, there is the identity
 \begin{equation}\label{TK1}
 H^x(t)= \lambda^x(t) + \sum_{0< s \leq t} \left|{\bf 1}_{[x,\infty)}V(s)- {\bf 1}_{[x,\infty)}V(s-)\right|,
 \end{equation}
 where $H^x(t)$ denotes the total variation on the time-interval $[0,t]$ of the step function ${\bf 1}_{[x,\infty)}\circ V$.

\end{document}